%% file: july.2002.tex
\documentclass[12pt]{article}

\newcommand{\oo}{\mathbf 0}

\usepackage{latexsym}
\usepackage{amssymb}
\usepackage{amsmath}
\usepackage{fullpage}
\usepackage{eepic}

\DeclareMathOperator{\Poi}{Poisson}

\begin{document}

\begin{center}
{\large An Illuminating Counterexample} \\ Michael Hardy
\end{center}

Suppose that $X_1,\dots,X_n$ are independent random variables
with a normal (or ``Gaussian'') distribution with expectation
$\mu$ and variance $\sigma^2$.  A statistician who has observed
the values of $X_1,\dots,X_n$ must guess the values of $\mu$
and $\sigma^2$.  Among the statistically naive, it is sometimes
asserted that
\begin{eqnarray*}
S^2=\frac{1}{n-1}\sum_{i=1}^n \left(X_i-\overline{X}\,\right)^2,
\end{eqnarray*}
where $\overline{X}=(X_1+\cdots+X_n)/n$,
is a better estimator of $\sigma^2$ than is
$$
T^2=\frac{1}{n}\sum_{i=1}^n \left(X_i-\overline{X}\,\right)^2,
$$
because $S^2$ is ``unbiased'' and $T^2$ is ``biased.''
That means $E(S^2)=\sigma^2\neq E(T^2)$, i.e., an
``unbiased estimator'' is a statistic whose expected value is
the quantity to be estimated.

The goodness of an estimator is sometimes measured by the smallness
of its ``mean squared error,'' defined as
$E\left(\left([\text{estimator}]
-[\text{quantity to be estimated}]\right)^2\right)$.  By that
criterion the biased estimator $T^2$ would be better than the
unbiased estimator $S^2$, since
$$
E((T^2-\sigma^2)^2)<E((S^2-\sigma^2)^2),
$$
but the difference is so slight that no one's statistical conscience
is horrified by anyone's preferring $S^2$ over $T^2$. Besides,
the smallness of the mean squared error as a criterion for
evaluating estimators is not necessarily sacred anyway.

A more damning example, well-known among statisticians,
is described in \cite[p.~168]{romano}.
We have $X\sim\Poi(\lambda)$, so that
$P(X=x)=\lambda^x e^{-\lambda}/x!$ for $x=0,1,2,\dots$, and
$P(X=0)^2=e^{-2\lambda}$ is to be estimated.  Any unbiased
estimator $\delta(X)$ satisfies
$$
E(\delta(X))=\sum_{x=0}^\infty\delta(x)
\frac{\lambda^x\, e^{-\lambda}}{x!}=e^{-2\lambda}
$$
uniformly in $\lambda\geq 0$.  Clearly the only such function is
$\delta(x)=(-1)^x$.  Thus, if it is observed that $X=200$, so that
it is astronomically implausible that $e^{-2\lambda}$ is anywhere
near 1, the desideratum of unbiasedness nonetheless requires us to
use $(-1)^{200}=1$ as our estimate of $e^{-2\lambda}$.
And if $X=3$ is observed, the situation is even more absurd:
we must use $(-1)^3=-1$ as an estimate of a quantity that we
know to be in the interval $(0,1)$.
A far better estimator of $e^{-2\lambda}$ is the biased
estimator $e^{-2X}$ (which is the answer given by the
well-known method of maximum likelihood).

Here is a different counterexample, which the visually
inclined may find even more horrifying.  A light source
is at an unknown location $\mu$ somewhere in the disk
$D=\{\,(x,y):x^2+y^2\leq 1\,\}$ in the Euclidean plane
(see Figure 1).
 \begin{figure}
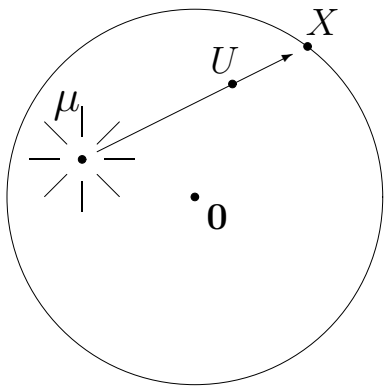

 \include{fig1}
 \caption{$D=\{\,(x,y):x^2+y^2=1\,\}$}
 \end{figure}
A dart thrown at the disk strikes some random
location $U$ in the disk, casting a shadow at a point
$X$ on the boundary.  The random variable $U$ is
uniformly distributed in the disk, i.e., the probability
that it is within any particular region is proportional
to the area of the region.  The boundary is a translucent
screen, so that an observer located outside of the disk
can see the location $X$ of the shadow, but cannot see
where either the light source or the opaque object is.
Given only that information---the location $X$ of the
shadow---the location $\mu$ of the light source must
be guessed.

A common-sense approach to guessing $\mu$ might proceed
as follows: {\em Before} we observe the shadow, our information
is invariant under rotations, and so should be our estimate.
Therefore, we use $0$ in $\mathbb{R}^2$ as our prior (i.e., pre-data)
estimate.  Then, when we observe $X$, since $X$ is more likely
to be far from the light source than close to it, we adjust our
estimate by moving it away from the shadow.  Because the amount
of information in the shadow is small, we don't move it very far.
We get an estimator of the form $cX$ with $c<0$, but $c$ is not
very much less than 0.

If we insist on unbiasedness, we must choose $c$ so that
$E(cX)=\mu$ uniformly in $\mu$.  To think about that, we
first express the problem in polar coordinates.
Write $\mu=\rho(\cos\varphi,\sin\varphi)$ and
$X=(\cos\Theta,\sin\Theta)$. \vspace{12pt} \\
{\bf Proposition:}
{\em The probability distribution
of the random angle $\Theta$ is given by}
\begin{equation}\label{conditional}
P(d\theta)
=\frac{1-\rho\cos(\theta-\varphi)}{2\pi} \, d\theta. 
\end{equation}

From this proposition it follows that $E(X)=-\mu/2$.
Therefore, our unbiased estimator is $cX=-2X$, which
is always absurdly remote from the $D$, by a full radius!
\vspace{12pt} \\
{\em Proof of the Proposition.}
A simplification will follow from the observation that the
way in which the probability distribution $P(d\theta)$
depends on $\mu$ is both rotation-equivariant and affine.
That it is affine means that if the probability distribution
of $\Theta$ is $P_{\mu}(d\theta)$ when the light source
is at $\mu$ then $P_{a\mu+(1-a)\nu}(d\theta)
=aP_{\mu}(d\theta)+(1-a)P_{\nu}(d\theta)$
for any value of $a$ for which $a\mu+(1-a)\nu$ remains
within the disk.  (An affine mapping is one that preserves
linear combinations in which the sum of the coefficients is 1;
a linear combination satisfying that constraint is a
``affine combination.'')  To see that this mapping is affine,
\begin{figure}
\include{fig2}
\caption{}
\end{figure}
consider Figure 2.  The area between $\mu$ and the arc from
$A$ to $B$ is the sum of the area of the triangle $\mu AB$
and the area of the region bounded by the arc $AB$ and the
secant line $AB$.  As $\mu$ moves, the area bounded by the
arc and the secant line remains constant and the area of the
triangle depends on $\mu$ in an affine fashion.  The desired
``affinity'' follows.

Rotation-equivariance reduces the problem to finding the
probability distribution when $\mu$ is between $(0,0)$ and
$(1,0)$.  ``Affinity'' reduces it from there to the problem
of finding the probability distribution when $\mu$ is at
either of those two points.

If $\mu=(0,0)$, the probability distribution of $\Theta$ is
clearly uniform on the interval from $0$ to $2\pi$, i.e., it
is $d\theta/(2\pi)$.
If $\mu=(1,0)$, then for $0\leq\theta\leq 2\pi$ we have
$$
P(0\leq\Theta\leq\theta)
=\frac{
\text{area between arc and straight line from }
(1,0)\text{ to }(\cos\theta,\sin\theta)}{\text{area of disk}}
=\frac{\theta-\sin\theta}{2\pi}.
$$
Differentiation yields
$$
P(d\theta)=\frac{1-\cos\theta}{2\pi}d\theta.
$$
If $\mu=(\rho,0)$ then by ``affinity'' we have
$$
P(d\theta)=(1-\rho)\frac{d\theta}{2\pi}
+\rho\frac{(1-\cos\theta)d\theta}{2\pi}
=\frac{1-\rho\cos\theta}{2\pi}d\theta.
$$
Rotation-equivariance then gives (\ref{conditional}).
$\quad\blacksquare$

The Bayesian approach to statistical inference assigns
probabilities, not to {\em events} that are {\em random},
according to their relative frequencies of occurrence,
but to {\em propositions} that are {\em uncertain}, according
to the degree to which known evidence supports them.
Accordingly, we can regard the location $\mu$ of the light
source as uniformly distributed in the disk, and then
use the conditional expected location $E(\mu\vert X)$ as
an estimator of $\mu$.  Equation (\ref{conditional})
gives the conditional distribution of $\Theta$ given $\mu$;
the marginal (i.e., ``unconditional'') distribution of
$\mu=\rho(\cos\varphi,\sin\varphi)$ is given by
\begin{equation}\label{marginal}
\frac{\rho\,d\rho\,d\varphi}{\pi}
\end{equation}
The joint distribution of $(\mu,\Theta)$ is the product
of (\ref{conditional}) and (\ref{marginal}):
\begin{equation}\label{joint}
\frac{(1-\rho\cos(\theta-\varphi))
\rho\,d\rho\,d\varphi\,d\theta}{2\pi^2}.
\end{equation}
The conditional distribution of $\mu=\rho(\cos\varphi,\sin\varphi)$
given that $\Theta=\theta$ comes from regarding (\ref{joint}) as a
function $\rho$ and $\varphi$ with $\theta$ fixed and normalizing:
$$
P(d\rho,d\varphi\vert\Theta=\theta)=\frac{
(1-\rho\cos(\theta-\varphi))\rho\,d\rho\,d\varphi}{
\text{constant}}.
$$
Integration shows that the ``constant'' is $\pi$.
Finally, we get
\begin{eqnarray*}
E(\mu \vert X) & = & \int_0^{2\pi} \int_0^1
\rho (\cos\varphi,\sin\varphi)
\frac{1-\rho\cos(\Theta-\varphi)}{\pi}\rho \, d\rho \, d\varphi \\
& = & -(\cos \Theta,\sin \Theta)/4 \ = \ -X/4,
\end{eqnarray*}
which is an eminently reasonable estimator under the
circumstances.

\noindent
{\em Department of Mathematics \\
Massachusetts Institute of Technology \\
Cambridge, MA 02139 \\
hardy@math.mit.edu}

\end{document}

%% file: fig1.tex
\setlength{\unitlength}{1mm}
\begin{picture}(162,30)(10,20)

\put(92,30){\circle{50}}
\put(95,29){\makebox(0,0)[t]{\large{$\oo$}}} 
\put(92,30){\circle*{1}}


\put(77,35){\circle*{1}}
\put(80,35){\line(1,0){4}}
\put(74,35){\line(-1,0){4}}
\put(77,38){\line(0,1){4}}
\put(77,32){\line(0,-1){4}}
\put(79,37){\line(1,1){3}}
\put(79,33){\line(1,-1){3}}
\put(75,37){\line(-1,1){3}}
\put(75,33){\line(-1,-1){3}}
\put(79,36){\vector(2,1){26}}
\put(75,44){\makebox(0,0)[t]{\Large{$\mu$}}}

\put(97,45){\circle*{1}}
\put(96,50){\makebox(0,0)[t]{\large{$U$}}}

\put(107,50){\circle*{1}}
\put(109,55){\makebox(0,0)[t]{\large{$X$}}}

\end{picture}

%% file: fig2.tex
\setlength{\unitlength}{1mm}
\begin{picture}(150,40)(-5,20)
\put(76,30){\circle{50}}
\put(76,45){\circle*{1}}
\put(76,45){\line(5,-3){25}}
\put(101,30){\line(-1,3){5}}
\put(96,45){\line(-1,0){20}}
\put(72,47){\makebox(0,0)[t]{\Large{$\mu$}}}
\put(104,31){\makebox(0,0)[t]{\large{$A$}}}
\put(100,47){\makebox(0,0)[t]{\large{$B$}}}

\end{picture}